\documentclass[11pt]{article}
\usepackage{amsfonts,amsmath,latexsym,color,epsfig,enumitem}
\date{}

\title{Disjointness graphs of short polygonal chains}

\author{
{\sl J\'anos Pach}\thanks{
R\'enyi Institute, Budapest and MIPT, Moscow; \texttt{pach@renyi.hu}; \texttt{pach@cims.nyu.edu}.
Supported by the National Research, Development and Innovation Office (NKFIH) grant K-131529, ERC Advanced Grant ``GeoScape,'' the Austrian Science Fund grant Z 342-N31, and by the Ministry of Education and Science of the Russian Federation in the framework of MegaGrant No.\ 075-15-2019-1926.}
\and
{\sl G\'abor Tardos}\thanks{
R\'enyi Institute, Budapest and MIPT, Moscow; \texttt{tardos@renyi.hu}.
Supported by the ERC Synergy Grant “Dynasnet”
No. 810115, the ERC advanced grant “GeoSpace” No. 882971, the National Research,
Development and Innovation Office — NKFIH projects K-132696, 
SSN-135643 and by the grant of Russian Government N 075-15-2019-1926.}
\and
{\sl G\'eza T\'oth}\thanks{
R\'enyi Institute, Budapest; \texttt{geza@renyi.hu}.
Supported by National Research, Development and Innovation Office, NKFIH, K-131529 and  ERC Advanced Grant ``GeoScape,''.}}

\date{}

\begin{document}

\maketitle

\begin{abstract}
The {\em disjointness graph} of a set system is a graph whose vertices are the sets, two being connected by an edge if and only if they are disjoint. It is known that the disjointness graph $G$ of any system of segments in the plane is {\em $\chi$-bounded}, that is, its chromatic number $\chi(G)$ is upper bounded by a function of its clique number $\omega(G)$.

Here we show that this statement does not remain true for systems of polygonal chains of length $2$. We also construct systems of polygonal chains of length $3$ such that their disjointness graphs have arbitrarily large girth and chromatic number. In the opposite direction, we show that the class of disjointness graphs of (possibly self-intersecting) \emph{$2$-way infinite} polygonal chains of length $3$ is $\chi$-bounded: for every such graph $G$, we have $\chi(G)\le(\omega(G))^3+\omega(G).$
\end{abstract}

\section{Introduction}
Ramsey theory has many applications to other parts of mathematics and computer science~\cite{R04}, including complexity theory~\cite{MeW87}, approximation algorithms, \cite{MS85}, coding~\cite{KSV05}, geometric data structures~\cite{MeN07}, graph drawing and representation~\cite{BEF98}. Constructing nearly optimal Ramsey graphs is a notoriously difficult combinatorial problem~\cite{F77}. The few efficient constructions that we have are far from optimal, but they can come in handy in those areas where we have interesting theorems, but lack nontrivial constructions. Here we provide two examples from combinatorial geometry, based on two classical constructions of Erd\H os and Hajnal~\cite{EH66, EH68}. We close this paper with a result pointing in the opposite  direction.

\smallskip

For any graph $G$, let $\chi(G)$ and $\omega(G)$ denote the {\em chromatic number} and the {\em clique number} of $G$, respectively. Clearly, we have $\chi(G)\ge\omega(G)$, and if equality holds for every induced subgraph of $G$, then $G$ is called a {\em perfect graph}.
Following Gy\'arf\'as and Lehel \cite{GL83}, \cite{GL85}, \cite{G85}, \cite{G87}, a class of graphs $\cal G$ is said to be {\em $\chi$-bounded} if there is a function $f$ such that $\chi(G)\le f(\omega(G))$ for every $G\in{\cal G}$.

Which classes of graphs are $\chi$-bounded? Or, reversing the question, if a graph has small clique number, how can its chromatic number be large? These questions are related to the some of the deepest unsolved problems in graph theory. There are two different approaches that have yielded spectacular results in recent years.
\smallskip

One can investigate what kind of substructures must necessarily occur in graphs of high chromatic number. According to Hadwiger's conjecture ~\cite{H43}, if the chromatic number of a graph is at least $t$, then it must contain a $K_t$-minor. (We now know that it contains a $K_s$-minor with $s>t/(\log\log t)^{18}$; cf.~\cite{Po20}.) Gy\'arf\'as~\cite{G75} proved that if a graph has bounded clique number and its chromatic number is sufficiently large, then it must contain a long induced path; see also~\cite{GrHM03}. According to the (still open) Gy\'arf\'as-Sumner conjecture \cite{Su81}, the same is true for any fixed tree instead of a path. Scott and Seymour proved that the class of graphs with no induced odd cycle of length at least $5$ is $\chi$-bounded. For many beautiful recent results of this kind, see the survey~\cite{SS20}.
\smallskip

The second fruitful research direction was initiated by Asplund and Gr\"unbaum~\cite{AG60}: Find geometrically defined classes of graphs that are $\chi$-bounded. Given a set $S$ of geometric objects, their {\em intersection graph (resp., disjointness graph)} is a graph on the vertex set $S$, in which two vertices are connected by an edge if and only if the corresponding objects have a nonempty intersection (resp., are disjoint). It was proved in~\cite{AG60} that the class of intersection graphs of axis-parallel rectangles in the plane is $\chi$-bounded (see also~\cite{ChW20}).  The corresponding statement is false for boxes in $3$ and higher dimensions~\cite{B65}, and even for segments in the plane~\cite{PKK14}.
\smallskip

For {\em disjointness} graphs $G$ of systems of segments in the plane, we have $\chi(G)\le (\omega(G))^4.$ The same is true for systems {\em $x$-monotone curves}, that is, for continuous curves in the plane with the property that every vertical line intersects them in at most one point. It was shown in~\cite{PaT20} that, in this generality, the order of magnitude of this bound cannot be improved. On the other hand, we proved~\cite{PaTT17} that the class of disjointness graphs of {\em strings} (continuous curves in the plane) is not $\chi$-bounded. Improving our construction, M\"utze, Walczak, and Wiechert~\cite{MWW18} exhibited systems of polygonal curves consisting of {\em three} segments such that their disjointness graphs are triangle-free ($\omega=2$), yet their chromatic numbers can be arbitrarily large.

The above results leave open the case of polygonal curves consisting of {\em two} segments. Our first result settles this case. A polygonal curve consisting of $k$ segments is called a {\em polygonal $k$-chain}.
\medskip

\noindent{\bf Theorem 1.} {\em There exist arrangements of polygonal $2$-chains in the plane whose disjointness graphs are triangle-free and have arbitrarily large chromatic numbers.}
\medskip

We do not know if Theorem 1 can be strengthened by requiring that the disjointness graph of the curves has large girth.
\medskip

\noindent{\bf Problem 2.} \emph{Do there exist arrangements of polygonal $2$-chains in the plane whose disjointness graphs have arbitrarily large girth and chromatic number?}
\medskip

Our next result shows that the answer to the above question is in the affirmative if, instead of $2$-chains, we are allowed to use polygonal $3$-chains.
\medskip

\noindent{\bf Theorem 3.} \emph{For any integers $g$ and $k$, there is an arrangement of non-selfintersecting polygonal $3$-chains in the plane whose disjointness graph has girth at least $g$ and chromatic number at least $k$.}
\medskip

A \emph{$1$-way infinite} polygonal $2$-chain is the union of a half-line and a segment that share an endpoint. In our proof of Theorem 1, we actually construct arrangements of $1$-way infinite polygonal $2$-chains whose disjointness graphs are triangle free, but have arbitrarily large chromatic number. Doubly tracing these $1$-way infinite $2$-chains and slightly perturbing the resulting curve, we obtain an arrangement of \emph{$2$-way infinite} $4$-chains, \emph{i.e.,} $4$-chains whose first and last pieces are half-lines. Hence, we obtain the following
\medskip

\noindent{\bf Corollary 4.} {\em There exist arrangements of $2$-way infinite polygonal $4$-chains in the plane whose disjointness graphs are triangle-free and have arbitrarily large chromatic numbers.}
\medskip

Our next theorem shows that Corollary 4 is optimal: the class of disjointness graphs of (possibly self-intersecting) $2$-way infinite polygonal $3$-chains is $\chi$-bounded.

\medskip
\noindent{\bf Theorem 5.} \emph{Let $G$ be the disjointness graph of an arrangement of $2$-way infinite polygonal $3$-chains in the plane.
Then we have $\chi(G)\le(\omega(G))^3+\omega(G)$.}
\medskip

In fact, we will establish Theorem 5 in a somewhat stronger setting: for arrangements of $2$-way infinite curves that consist of three $x$-monotone pieces; see Theorem~7. With more work, the bound in Theorem~5 and Theorem~7 can be improved to $\chi(G)\le(\omega(G))^3$.

In the
polygonal case, our proof is algorithmic. There is a polynomial time
algorithm in the number of the polygonal chains, which, for every $k$,
either finds $k$ pairwise disjoint chains or produces a coloring of
their disjointness graph with at most $k^3$ colors.

\smallskip

In Sections~2 and~3, we establish Theorems~1 and~3, respectively. Section~4 contains the proof of Theorem~5. We end this note with a few remarks and open problems.
\smallskip

In the sequel, we informally call a polygonal $2$-chain a  \emph{V-shape} and a polygonal $3$-chain a \emph{Z-shape}.

\section{Shift graphs---Proof of Theorem 1}

For every $n>1$, Erd\H os and Hajnal~\cite{EH68} defined the {\em shift graph} $S_n$, as follows.
The vertex set of $S_n$ consist of all pairs $(a,b)$ with $1\le a<b\le n$, where
two vertices, $(a, b)$ and $(a',b')$, are connected by an edge if and only if
$b=a'$ or $b'=a$. It is easy to see that $S_n$ is triangle-free and that $\chi(S_n)=\lceil\log_2n\rceil$.
\smallskip

Order the vertices $(a,b)$ of $S_n$ according to the {\em co-lexicographic order}, that is,
let $(a, b)\prec (a', b')$ if $b<b'$, or if $b=b'$ and $a<a'$.
Let $v_1, \ldots, v_{\binom{n}{2}}$ denote the vertices of $S_n$, listed in this order.

Let $v_i=(a, b)$ be a vertex. Its neighbors having a smaller
index are $(a',b')$ with $b'=a$. No such neighbor exist if and only if $a=1$.
Notice that, for any $i$,
\begin{enumerate}
\item {\em either} $v_i$ has no neighbor $v_j$ with a smaller index $j<i$,
\item {\em or} there exist integers $c(i), d(i)$ with $1\le c(i)\le d(i)<i$ such that
for every $j<i$, $$v_jv_i\in E(S_n)\Longleftrightarrow c(i)\le j\le d(i).$$
\end{enumerate}

Recall that a {\em $1$-way infinite V-shape} is the union of a half-line and a segment that share an
endpoint. In the rest of this proof, for simplicity, we call a $1$-way infinite V-shape \emph{long}.

Our goal is to assign a long V-shape to each vertex of $S_n$ so that two
V-shapes are disjoint if and only if the corresponding vertices are adjacent in $S_n$. This will prove Theorem~1, because in any finite collection of long V-shapes, we can cut the half-lines short so that the resulting (bounded) V-shapes have the same intersection structure. Hence, we obtain a collection of V-shapes with $S_n$ as its disjointness graph, and the graphs $S_n$ are triangle-free and their chromatic numbers tend to infinity, as $n\rightarrow\infty$.
\smallskip

We assign the long V-shape $V_i$ to the vertex $v_i$ of $S_n$ recursively starting at $V_1$.
Let $h_i$ and $s_i$ denote the half-line and the straight-line segment, respectively,
comprising $V_i$ and let us denote their common endpoint by $p_i=(x_i, y_i)$.
We write $q_i$ for the other endpoint of $s_i$.

During the recursive process, we will maintain the following properties:
\smallskip

\emph{(i)} $p_i$ is the left end point of both $h_i$ and $s_i$;
\smallskip

\emph{(ii)} both $h_i$ and $s_i$ have positive slopes;
\smallskip

\emph{(iii)} $s_i$ is above $h_i$, {\em i.e.},
the slope of $s_i$ is larger than the slope of $h_i$;
\smallskip

\emph{(iv)} for any $i>j$, the slope of $h_i$ will be smaller than the slope of $h_j$;
\smallskip

\emph{(v)} for any $i>j$, we have $x_i<x_j$ and $y_i<y_j$.

\smallskip

Let $V_1$ be any long V-shape satisfying the above conditions.
Let $i>1$, and assume recursively that we have already constructed the long V-shapes
$V_1,\ldots, V_{i-1}$ satisfying the above requirements. Next, we define $V_i$.
We distinguish two case;
\smallskip

{\bf Case A:} {\em The vertex $v_{i}=(a,b)$ has no neighbor with a smaller index, {\em i.e.,} we have $a=1$.}

Let $\ell$ be a horizontal line passing above $p_1$. It will intersect every $V_j$ with $1\le j< i$.
Slightly rotate $\ell$ about any fixed point of the plane
so that the resulting line $\ell'$ has a very small positive slope, smaller than the slope of $h_{i-1}$ and it still intersects all $V_j$ for $j<i$.
Choose a point $p_{i}=(x_i,y_i)\in \ell,$, very far to the left, so that $x_i<x_{i-1}$ and $y_{i}<y_{i-1}$.
Let $h_{i}$ be the part of $\ell'$ to the right of $p_{i}$,
and let $q_{i}$ be a point to the right of $p_{i}$ which lies above $h_{i}$. One can choose $q_i$ such that the segment $s_i=p_iq_i$ does not intersect any of the earlier $V_j$.
\smallskip


\begin{figure}[!ht]
\begin{center}
\scalebox{0.3}{\includegraphics{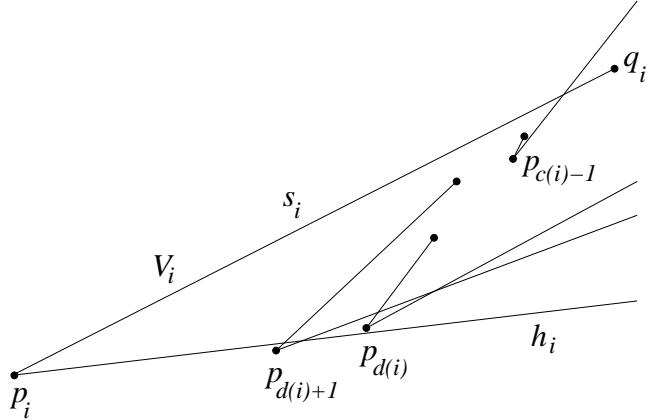}}
\caption{Inserting $V_i$.}\label{vbetu}
\end{center}
\end{figure}

{\bf Case B:} {\em The vertex $v_{i}=(a,b)$ has at least one neighbor of smaller index, {\em i.e.,} $a>1$.}

Let $c(i)$ and $d(i)$ be the constants satisfying property~(2) above and let $\ell$ be a horizontal line that passes below $p_{d(i)}$ and
above $p_{d(i)+1}$. In case $d(i)+1=i$ we could simply choose $\ell$ to be an arbitrary horizontal line below $p_{d(i)}$, but the careful reader may notice that this case never occurs as no vertex $v_i$ in $S_n$ is adjacent to $v_{i-1}$.

The line $\ell$ intersects every $V_j$ with $d(i)<j<i$ and is disjoint from all $V_j$ with $j\le d(i)$.
Slightly rotate $\ell$ about any fixed point in the plane so that the resulting line $\ell'$
has a very small positive slope, smaller than that of $h_{i-1}$ and it still intersects the same previously defined long V-shapes $V_j$. Select a slope $\alpha$ which is larger than the slope of $h_{c(i)}$,
but smaller  than the slope of $h_{c(i)-1}$, if $h_{c(i)-1}$ exists, that is,
if $c(i)>1$.

For any $j<i$, let $\ell_j$ and $\ell'_j$
denote the lines of slope $\alpha$ through $p_j$ and $q_j$, respectively.
Choose a point $p_{i}=(x_i,y_i)\in\ell'$ so far to the left
that we have $x_i<x_{i-1}$, $y_{i}<y_{i-1}$ and $p_i$ lies above the lines $\ell_j$ and $\ell'_j$, for all $j\le i$.

Let $h_{i}$ be the part of $\ell$ to the right of $p_{i}$.
Let $f$ be the half-line of slope $\alpha$, whose left endpoint is $p_{i}$.
Then $f$ goes strictly above all $s_j$ for $j<i$ and also of all $h_j$ with $c(i)\le j<i$, but will intersect all $h_j$ with $1\le j<c(i)$.
Choose $q_{i}$ on $f$ to the right of these intersection points, then the segment $s_i=p_iq_i$ also intersects all $h_j$ with $1\le j<c(i)$.

Notice that the long V-shape $V_i$ consisting of $h_i$ and $s_i$ constructed above satisfies the conditions (i)--(v) listed above, further it intersects exactly those other long V-shapes $V_j$ ($j<i$) for which $v_j$ and $v_i$ are not adjacent in $S_n$.
See Fig. 1.
This means that the disjointness graph of the collection of the $\binom n2$ long V-shapes constructed above is exactly $S_n$. This completes the proof of Theorem~1. \hfill $\Box$
\medskip

In the above proof, we have constructed a collection of $1$-way infinite V-shapes in which each pair intersects at most twice. With a little additional care (namely, by insisting that each $q_i$ is higher than $p_1$), we can achieve the following. For $1\le i<j\le\binom n2$, with $v_i=(a,b)$ and $v_j=(a',b')$, we have
\begin{itemize}
\item if $a'<b$, then $V_i$ and $V_j$ intersect once;
\item if $a'=b$, then $V_i$ and $V_j$ are disjoint;
\item if $a'>b$, then $V_i$ and $V_j$ intersect twice.
\end{itemize}

\section{Hypergraphs of large girth---Proof of Theorem~3}

A hypergraph $H$ is a pair $(V,E)$, where $V$ is a finite vertex set,
$E$ is the set of hyperedges, that is, a collection of subsets of $V$.
$H$ is called \emph{$n$-uniform} if each of its hyperedges has $n$ vertices.
In a \emph{proper coloring} of $H$, every vertex is  assigned
a color in such a way that none of the hyperedges is monochromatic.
The \emph{chromatic number} of $H$ is the smallest number of colors
used in a proper coloring of $H$.
A \emph{Berge-cycle} in $H$ consists of a sequence of
distinct vertices $v_1,\ldots,v_k$ and a sequence of distinct
hyperedges $e_1,\ldots,e_k\in E$ with $v_i,v_{i+1}\in e_i$ for
$1\le i<k$ and $v_k,v_1\in e_k$. Here $k$ is the \emph{length}
of the Berge-cycle and it is assumed to be at least $2$.
The \emph{girth} of a hypergraph is the length of
its shortest Berge-cycle (or infinite if it has no Berge-cycle).

For the proof, we need the following classical result.


\medskip

\noindent{\bf Erd\H os-Hajnal Theorem}  \cite{EH66} (Corollary 13.4).
\emph{For any integers $n\ge2$, $g\ge3$, and $k\ge2$, there exists
  an $n$-uniform hypergraph with girth at least $g$ and chromatic number at least $k$.}
\medskip

Theorem 3 is a direct consequence of part (5) of the following statement.
\medskip

\noindent {\bf Lemma 6.}  {\em For any integers $g\ge3$, $k\ge2$, there is a natural number $n=n(g,k)$ such that
for every set $P$ of $n$ points on the $x$-axis in $\bf{R}^2$
and for every real $c>0$, there is an arrangement $Z=Z(P)$ of $n$ Z-shapes satisfying the following conditions.
\smallskip

\emph{(1)} Each point in $P$ is the endpoint of exactly one Z-shape in $Z$.
\smallskip

\emph{(2)} Apart from a single endpoint in $P$, every Z-shape in $Z$ lies strictly above the $x$-axis.
\smallskip

\emph{(3)} No Z-shape in $Z$ is self-intersecting and any two cross at most twice.
\smallskip

\emph{(4)} 
For any Z-shape $pqrs\in Z$ with $p\in P$, whose vertices have $x$-coordinates $x_p$, $x_q$, $x_r$, and $x_s$, respectively, we have $x_q+c<x_p<x_s<x_r-c$.

\smallskip

\emph{(5)} The disjointness graph of the Z-shapes
in $Z$ has girth at least $g$ and chromatic number at least $k$.}
\medskip

\textbf{Proof.}
For each $g$, we prove the lemma by induction on $k$.
We fix $g\ge 3$. For $k=2$,  $n(g,2)=2$ is a good choice. For any two points on the $x$-axis and any $c>0$, we can take
two disjoint Z-shapes satisfying the requirements. Their disjointness graph is $K_2$, its chromatic number $2$ and
it has infinite girth. See Fig. 2.


\begin{figure}[!ht]
\begin{center}
\scalebox{0.4}{\includegraphics{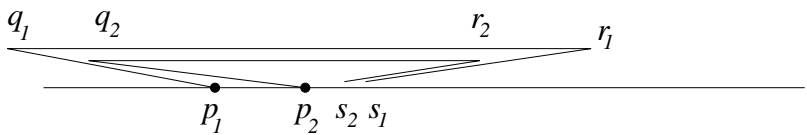}}
\caption{The case $k=2$.}\label{zbetu2}
\end{center}
\end{figure}

Suppose now that $k\ge 2$ and that we have already proved the statement for $k$. Now we prove it for $k+1$.
Let $n=n(g,k)$.

By the Erd\H os-Hajnal Theorem stated above, there exists an $n$-uniform hypergraph $H$ whose
girth is at least $g$ and chromatic number at least $k+1$.
Let $v_1, v_2,\dots, v_m$ denote the vertices of $H$ and $e_1, e_2, \ldots, e_M$ the hyperedges of $H$. Let $N=nM+m$. We show that $n(g,k+1)=N$
satisfies the requirements of the lemma.
\medskip

Let $P$ be an arbitrary
set of $N$ points on the $x$-axis and let $c>0$.
For any $v_i\in V(H)$, let $d_i$ denote the
\emph{degree} of $v_i$, that is, the number of hyperedges that contain $v_i$.
Obviously, we have
$$\sum_{i=1}^m(d_i+1)=nM+m=N.$$

Choose $m$ disjoint open intervals, $I_1, \ldots, I_m,$ such that each $I_i$ contains precisely $d_i+1$ points of $P$.
For every $i$, $1\le i\le m$, we associate the interval $I_i$ with vertex $v_i$ of $H$.
Let $p_i$ denote the leftmost point in $P\cap I_i$.
For every $i$ and  $j$ ($1\le i\le m$, $1\le j\le M$) for which $v_i\in e_j$,
assign a distinct point  $p_i^j \in (P\cap I_i)\setminus\{p_i\}$ to the pair $(v_i,e_j)$.
\medskip

Next, we construct a set of $N$ Z-shapes that satisfy conditions (1)--(5) of the lemma with parameters $g, k+1$, and $c$.
We construct subsets $Z_j$ of our eventual set of Z-shapes for $1\le j\le M$. We construct these sets one by one starting at $Z_1$ and using the inductive hypothesis for various subsets of $P$ of size $n$ and with a parameter $c'$ that we choose to be larger than $c$ plus the diameter of $P$.

For $j=1$, consider the  $n=n(g,k)$-element point set $P'_1=\{p_i^1 : v_i\in e_1\}$. By the induction hypothesis,
there is a set $Z_1$ of Z-shapes such that one of their endpoints belongs to $P'$, and they satisfy conditions~(1)--(5) with parameter $c'$.

Suppose that $j>1$ and that we have already constructed the sets of Z-shapes $Z_1, \ldots, Z_{j-1}$.
Now let $P'_j=\{p_i^j : v_i\in e_j\}$. By the induction hypothesis, there is a set $Z'_j$ of
Z-shapes with one of their endpoints in $P'$ which satisfy conditions~(1)--(5) with parameter $c'$.
Apply an affine transformation $(x,y)\to(x, y/K_j)$ to all Z-shapes in $Z'_j$, where $K_j$ is a very large constant to be specified later.
The resulting family of Z-shapes, $Z_j$, still satisfies all defining conditions and,
by choosing $K_j$ large enough, we can achieve that every element of $Z_j$ intersects
every Z-shape in $\bigcup_{h<j}Z_h$ exactly once or twice.
\smallskip

The set $\bigcup_{j=1}^MZ_j$ contains exactly one Z-shape starting at each point $p_i^j$.
We still need to add one Z-shape $z_i=p_iq_ir_is_i$ starting at each point $p_i,\; 1\le i\le m$.
We define them recursively for $i=1,\ldots,m$. We make sure that each $z_i=p_iq_ir_is_i$ satisfies the following properties.
\medskip

(i)\,\, \emph{The segment $q_ir_i$ is horizontal and the $y$-coordinate of its points is larger than the
$y$-coordinate of any point of any Z-shape in $(\bigcup_{j=1}^{M}Z_j)\cup\{z_h : 1\le h<i\}$.}
\smallskip

(ii)\, \emph{The slope of $p_iq_i$ is $-\varepsilon_i$, the slope of $r_is_i$ is $\varepsilon_i$, for a sufficiently small
constant $\varepsilon_i>0$, to be specified later.}
\smallskip

(iii) \emph{The $x$-coordinate of $s_i$ is equal to the $x$-coordinate of the right endpoint of $I_i$,
and the $y$-coordinate of $s_i$ is $\varepsilon_i$.}

\medskip

Clearly, if we choose $\varepsilon_i>0$ sufficiently small,
then $z_i$ is disjoint from all Z-shapes in $\bigcup_{j=1}^MZ_j$ that start
in $I_i$, but it intersects exactly once all other Z-shapes already defined.
Also, $z_i$ satisfies conditions~(2) and (3), and it satisfies condition (4), too,
provided that $\varepsilon_i$ is sufficiently small. See Fig. 3.


\begin{figure}[!ht]
\begin{center}
\scalebox{0.4}{\includegraphics{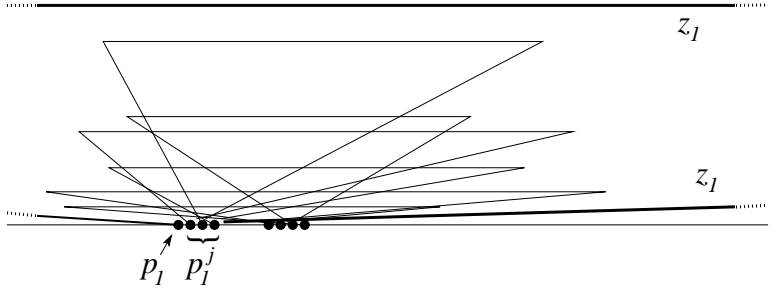}}
\caption{Inserting $z_1$.}\label{zbetu}
\end{center}
\end{figure}

As we maintained conditions~(1)--(4) throughout the construction, it remains only to prove
that the disjointness graph $G$ of $Z$ satisfies condition~(5) with $k+1$ in place of $k$.

To this end, let us explore the structure of $G$.
The vertices of $G$ can be partitioned into the sets $Z_j$ for $1\le j\le M$
and the independent set $W=\{z_i : 1\le i\le m\}$.
Further, there is no edge between two distinct sets $Z_j$ and $Z_{j'}$.
There is a single edge from $z_i$ to $Z_j$ if $v_i\in e_j$, and there is no edge
from $z_i$ to $Z_j$ otherwise.
Finally, each vertex in $Z_j$ is adjacent to exactly one of the vertices
$z_i$, and it satisfies $v_i\in e_j$.
\smallskip

The structure above implies that each cycle $C$ of $G$ is either contained in a
single set $Z_j$, or it passes through several sets $Z_j$ and several
vertices in $W$. In the former case,
by our assumption on the disjointness graph of $Z_j$, the length of $C$ is at most $g$.
In the latter case, let us record the vertices of $W$ and the sets $Z_j$
as the cycle passes through them: $z_{i_1},Z_{j_1},z_{i_2},Z_{j_2},\dots,z_{i_h},Z_{j_h}$.
Here, the vertices $v_{i_1},\dots,v_{i_h}$ are all distinct
and, if the same is true for the hyperedges $e_{j_1},\dots,e_{j_h}$,
then they form a Berge-cycle of length $h$ in the hypergraph $H$.
If the hyperedges are not all distinct, then an even shorter
Berge-cycle is formed by any repetition-free interval between
two occurrences of the same hyperedge.
By our assumption on the girth of $H$, we have $h\ge g$ in both cases,
so all cycles of $G$ have length at least $g$, as required.
\smallskip

Suppose now that there is a proper $k$-coloring of $G$.
Restricting it to the set $W$ (and identifying each $z_i\in W$ with the vertex $v_i$ of $H$),
we obtain a $k$-coloring of the vertices of the hypergraph $H$.
By our assumption, this cannot be a proper coloring. Therefore, there is a monochromatic hyperedge $e_j$.
In this case, no vertex in $Z_j$ can receive the common color of the vertices of $e_j$,
so we have a proper $(k-1)$-coloring of $Z_j$.
This contradicts our assumption on the disjointness graph of $Z_j$ and,
thus, proves that $G$ has no proper $k$-coloring. This concludes the proof of Lemma~6 and, hence, of Theorem~3. \hfill $\Box$

\medskip

James Davies \cite{D21} used a very similar construction to show that there are
intersection graphs of axis-parallel boxes and intersection graphs of lines
in $3$-space with arbitrarily large girths and chromatic numbers.

\section{Two-way infinite polygonal chains---Proof of Theorem~5}

As we pointed out at the end of Section~2, the class of disjointness graphs of $1$-way infinite V-shapes is not $\chi$-bounded. But if we require both ends of a V-shape to be long, the situation will change.
\smallskip

A \emph{$2$-way infinite polygonal $k$-chain} is a continuous curve in the plane consisting of two half-lines connected by an
(ordinary) polygonal $(k-2)$-chain. We can relax this definition by requiring only that each of the $k$ pieces are $x$-monotone, and the first and the last pieces have unbounded projections to the $x$-axis. In this case, the curve is called a $2$-way infinite \emph{$k$-monotone chain}.

According to this definition, a $2$-way infinite polygonal $2$-chain (V-shape)) whose half-lines are not vertical is a $2$-way infinite \emph{$2$-monotone} chain. It can also be regarded as a degenerate $2$-way infinite \emph{$3$-monotone} chain. Note that by performing a suitable rotation, if necessary, we can always assume that none of the half-line pieces of a finite arrangement of $2$-way infinite polygonal $k$-chains is vertical. Therefore, the following theorem implies Theorem~5.

\vskip0.5cm

\begin{figure}[!ht]
\begin{center}
\scalebox{0.3}{\includegraphics{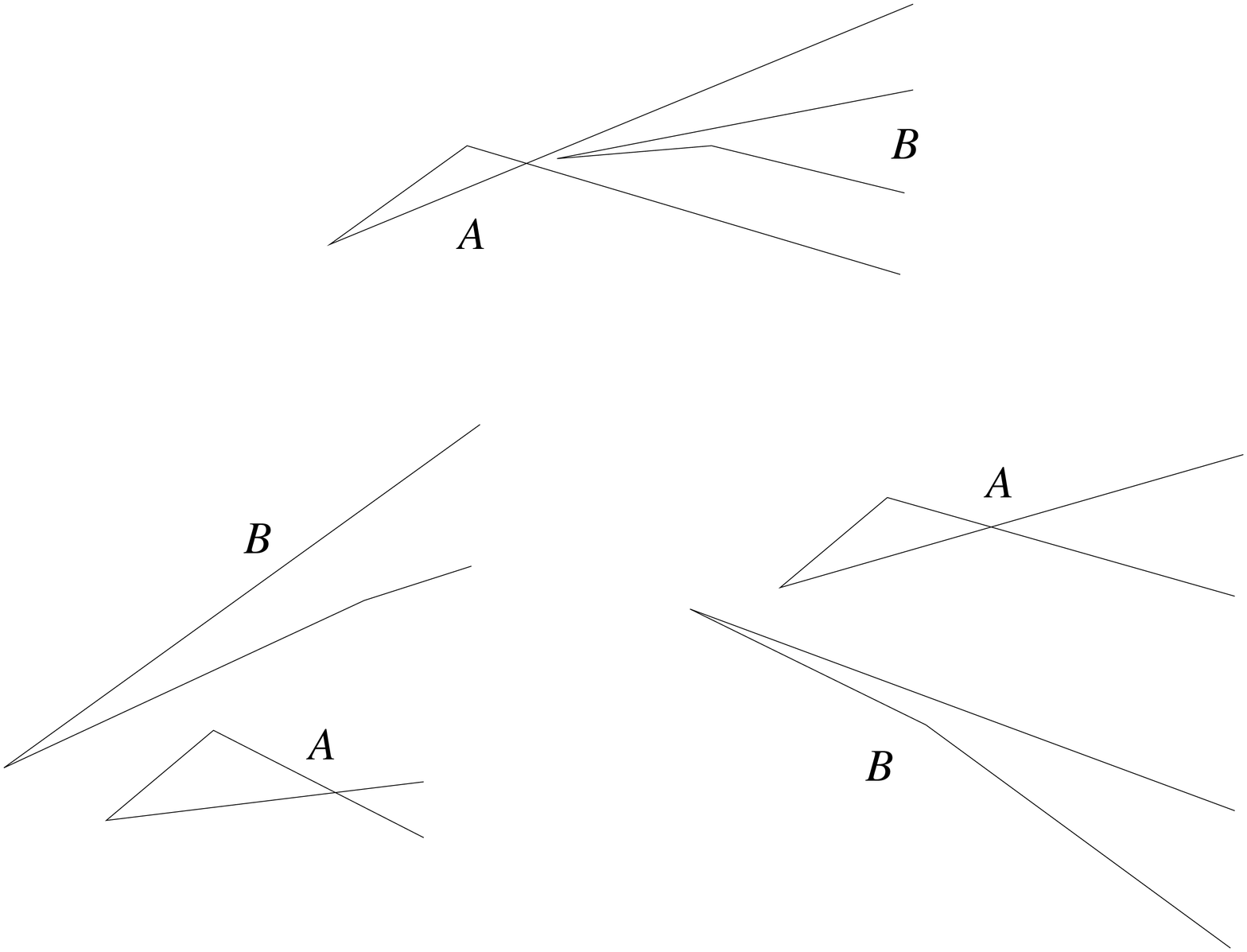}}
\caption{The three partial orders: $A$ is to the left of $B$, below $B$, and above $B$.}\label{3chain}
\end{center}
\end{figure}

\medskip

\noindent\textbf{Theorem 7.}  \emph{The disjointness graph $G$ of a finite arrangement of $2$-way infinite $3$-monotone chains satisfies $\chi(G)\le(\omega(G))^3+\omega(G)$.}
\medskip

\textbf{Proof.}
We call a (possibly self-intersecting) 2-way infinite $k$-monotone chain $A$ \emph{wide} if it intersects every vertical line. A chain $A$ with this property divides the plane into (open) connected components, exactly one of which contains a vertical half-line pointing upwards. We call this component the \emph{upside of $A$}. For any two wide 2-way infinite $k$-monotone chains, $A$ and $B$, we say that $A$ is \emph{higher} than $B$ if $A$ is contained in the upside of $A$. In this case, the upside of $B$ is also contained in the upside of $A$. Therefore, the relation ``higher'' defines a partial order on any arrangement of wide $k$-monotone chains. According to this partial order, only \emph{disjoint} pairs are comparable. Since any two disjoint \emph{wide} 2-way infinite $k$-monotone chains are comparable, the disjointness graph of any collection of wide 2-way infinite $k$-monotone chains is a comparability graph.
\smallskip

Now we turn our attention to the non-wide case. The complement of a non-wide 2-way infinite $k$-monotone chain $A$ has precisely one connected component which contains a vertical line. We call this component the \emph{large component}. The chain $A$ is said to be a \emph{right chain} if $A$ is to the right of the vertical lines in the large component, otherwise it is a \emph{left chain}. If $A$ is a right chain, we call its large component the \emph{left side} of $A$. On the other hand, if $A$ is a left chain, we call the union of \emph{all} connected components of the complement of $A$, other than its large component, the \emph{left side} of $A$.
\smallskip

For any two non-wide 2-way infinite $k$-monotone chains, $A$ and $B$, we say that $A$ is \emph{to the left} of $B$ if both $A$ and its left side are contained in the left side of $B$. Obviously, this relation also defines a partial order, with respect to which only disjoint non-wide chains are comparable. It is not true that any two disjoint non-wide 2-way infinite $3$-monotone chains are comparable. Therefore, we need to introduce two further partial orders.
\smallskip

For any two subsets of the plane, $A$ and $B$, we say that $A$ is
{\em below} $B$ ($A$ is {\em above} $B$, resp.), if the following two conditions are satisfied:
\begin{enumerate}
\item every vertical line that intersects $A$ also intersects $B$;
\item if $a\in A\cap\ell$ and $b\in B\cap\ell$ for a vertical line $\ell$, then the $y$-coordinate of $a$ is strictly \emph{lower} (\emph{higher}, resp.) than the $y$-coordinate of $b$.
\end{enumerate}
Note that ``above'' and ``below'' are two separate partial orders and not
the inverses of each other. 
It is clear that both of these relations are partial orders on arbitrary planar sets and that any two comparable sets are disjoint. See Fig. 4. 
\medskip

\noindent\textbf{Lemma 8.}  \emph{Any two disjoint non-wide 2-way infinite $3$-monotone chains, $A$ and $B$, are comparable by one of the three relations ``below'', ``above'', or ``to the left''.}
\medskip

To establish the lemma, note that non-wide 2-way infinite $3$-monotone chains must be, in fact, 2-way infinite \emph{$2$-monotone} chains. A left chain with this property is the union of the graphs of two continuous functions $f_1,f_2:(-\infty,a]\to\mathbb R$, where $f_1(a)=f_2(a)$. Let $B$ be another left chain obtained as the union of the graphs of two continuous functions $g_1,g_2:(-\infty,b]\to\mathbb R$, and assume that $A$ and $B$ are disjoint. We can assume, by symmetry, that $b\le a$. Consider $g_1(b)=g_2(b)$. It is easy to see that if it is below both $f_1(b)$ and $f_2(b)$, then $B$ is below $A$. If it is above both $f_1(b)$ and $f_2(b)$, then $B$ is above $A$. Finally, if $g_1(b)$ is between $f_1(b)$ and $f_2(b)$, then $B$ is to the left of $A$.  A similar argument applies if both $A$ and $B$ are right chains. Finally, if a left chain is disjoint from a right chain, then the left chain is always to the left of the right chain. This completes the proof of Lemma~8.
\smallskip

Now we return to the proof of Theorem~7. Fix a family $F$ of 2-way infinite $3$-monotone chains, and let $G$ denote their disjointness graph. Let $F_1$ and $F_2$ consist of the wide and non-wide elements of $F$, respectively. We have seen that the disjointness graph $G[F_1]$ of $F_1$ is a comparability graph. Comparability graphs are perfect, so we have $\chi(G[F_1])=\omega(G[F_1])$. We also proved that the comparability graph $G[F_2]$ of $F_2$ is the union of three comparability graphs. This implies that $\chi(G[F_2])\le(\omega(G[F_2]))^3$.

For the entire graph $G$, we have $$\chi(G)\le\chi(G[F_1])+\chi(G[F_2])\le\omega(G[F_1])+(\omega(G[F_2]))^3\le\omega(G)+(\omega(G))^3,$$ as required. This completes the proof of the theorem. \hfill $\Box$
\medskip

In \cite{PaT20}, for every $k\ge 2$, we constructed arrangements of $x$-monotone curves such that their left endpoints lie on the $y$-axis and their disjointness graphs have clique number $k$ and chromatic number ${k + 1\choose 2}.$ We can extend these curves to the left by adding horizontal half-lines without changing their intersection structure. Traversing the resulting curves twice, we obtain families of $2$-way infinite $2$-monotone chains such that their disjointness graphs satisfy $\chi(G)={\omega(G)+1\choose 2}$. We do not know whether the order of magnitude of the bound in Theorem~7 is best possible.

We were unable to improve on the bound in Theorem 7 even for 2-way infinite polygonal 3-chains. The best lower bound we have in this case is $ \omega(G)^{(\log 5/\log 2)-1}\approx \omega(G)^{1.32}$, and it follows from a construction in \cite{LMPT94}.

\section{Concluding remarks.}

{\bf A.} In Problem~2, we asked whether the disjointness graph of an arrangement of V-shapes can have simultaneously arbitrarily high chromatic number and girth. The following statement provides an affirmative answer to a relaxed version of this question. The \emph{odd-girth} of a graph is the length of the shortest odd cycle in it (or infinite if the graph is bipartite).
\medskip

\noindent\textbf{Proposition 9.} \emph{There exist arrangements of polygonal 2-chains in the plane whose disjointness
graphs have arbitrarily large odd-girths and chromatic numbers.}

\medskip

\textbf{Proof.} The proof is based on the same idea as the Proof of Theorem~1, where we represented the shift graph $S_n$ as the disjointness graph of an arrangement of V-shapes. The vertices of $S_n$ are pairs $(a,b)$ of integers $1\le a<b\le n$, so they can be associated with the edges of the complete graph $K_n$. Thus, the vertices of $S_n$ associated with the edges of a subgraph $G\subseteq K_n$ induce a subgraph $G^*\subseteq S_n$. It is easy to verify that for any $G\subseteq K_n$, we have
\smallskip

(1)  $\chi(G^*)\ge\log(\chi(G))$ \emph{and}
\smallskip

(2)  \emph{the odd-girth of $G^*$ is strictly larger than the odd-girth of $G$.}
\smallskip

For any integers $g$ and $k$, there exist $n=n(g,k)$ and a subgraph $G\subset K_n$ with girth (and, hence, odd-girth) at least $g$ and chromatic number at least $k$. By properties (1) and (2), the odd-girth of the corresponding induced subgraph $G^*$ of $S_n$ will be larger than $g$, and its chromatic number will be at least $\log k$. The graph $G^*$ inherits from $S_n$ a representation as a disjointness graph of V-shapes. \hfill $\Box$
\smallskip

Unfortunately, getting rid of short \emph{even} cycles, even 4-cycles, looks impossible by using this simple trick.
\medskip

{\bf B.} In the example of Theorem~1 and in all interesting known examples, whenever a class $\cal G$ of graphs,
closed under taking induced subgraphs is not $\chi$-bounded, then there exist {\em triangle-free} graphs $G\in{\cal G}$ with arbitrary large chromatic numbers. Is this true in general? Is it true that, for any $k$ and $\omega$, there exists $\chi=\chi(n,\omega)$ such that every graph $G$ with $\chi(G)\ge \chi$ and $\omega(G)\le\omega$ has an induced triangle-free subgraph $G'\subseteq G$ with $\chi(G')\ge k$? The existence of a not necessarily induced triangle-free subgraph $G'\subseteq G$ with this property was proved by R\"odl~\cite{R77}.
\medskip

{\bf C.} The arrangements of polygonal curves proving Theorems~1 and 3 have the property that any two of them have at most {\em two} points in common. It would be interesting to decide whether these theorems remain true if we insist that the curves are \emph{single-crossing}, that is, any two curves have at most one point in common at which they properly cross.
\medskip

\noindent\textbf{Conjecture 10.} \emph{The class of disjointness graphs of single-crossing polygonal $2$-chains is $\chi$-bounded.}
\smallskip

M\"utze \emph{et al.} \cite{MWW18} proved that the same statement is false for polygonal $3$-chains.
\medskip

{\bf D.} To prove Theorem~1, we established that the shift graph $S_n$, a triangle-free graph of unbounded chromatic number, can be obtained as the disjointness graph of V-shapes. However, the \emph{fractional chromatic number} of $S_n$ is bounded: it is smaller than $4$ for every $n$. Do there exist triangle-free disjointness graphs of V-shapes with arbitrarily large fractional chromatic number?

Analogously, our construction for Theorem 3 gives disjointness graphs with bounded fractional chromatic number. Do there exist disjointness graphs of Z-shapes with arbitrarily large girth and fractional chromatic number?
\medskip



\end{document}